\documentclass[hidelinks,onefignum,onetabnum]{siamart220329}

\usepackage{amsfonts,upquote}
\usepackage{graphicx,verbatim}

\headers{Polynomial and rational convergence rates}{Trefethen}

\title{Polynomial and rational convergence rates for Laplace problems on planar domains
\thanks{Submitted to the editors DATE.}}

\author{Lloyd N. Trefethen\thanks{School of Engineering and Applied Sciences,
Harvard University, Cambridge, MA 02138, USA}
(\email{trefethen@seas.harvard.edu}).}

\def\pput(#1,#2)#3{\noindent\smash{\raise#2pt\hbox to 0pt
   {\kern #1pt #3\hss}}\ignorespaces}

\def\complex{{\mathbb{C}}}
\def\En{E_n}
\def\Enn{E_{nn}}
\def\Uin{U_{\hbox{\scriptsize\rm in}}}
\def\Uout{U_{\hbox{\scriptsize\rm out}}}
\def\Re{\hbox{\rm Re\kern 1pt}}
\def\Im{\hbox{\rm Im\kern 1pt}}

\def\Oc{\overline{\Omega}}

\def\half{\textstyle{1\over 2}}
\def\pdes{PDE\kern .4pt s}
\def\ls{\llap{$*$\kern 1.5pt}}
\def\qed{~\hbox{\vrule width 3pt depth 4.5 pt height 2.5 pt}}
\def\Sc{\overline{S}\kern .7pt}

\begin{document}

\maketitle


\begin{abstract}
Laplace problems on planar domains can be solved by means of
least-squares expansions associated with polynomial or rational
approximations.  Here it is shown that, even in the context of
an analytic domain with analytic boundary data, the difference
in convergence rates may be huge when the domain is nonconvex.
Our proofs combine the theory of the Schwarz function for analytic
continuation, potential theory for polynomial and rational
approximation rates, and the theory of crowding of conformal maps.

\end{abstract}



\section{\label{intro}Introduction}
\begin{figure}[h!]
\begin{center}
\vskip 13pt
\includegraphics[scale=.85]{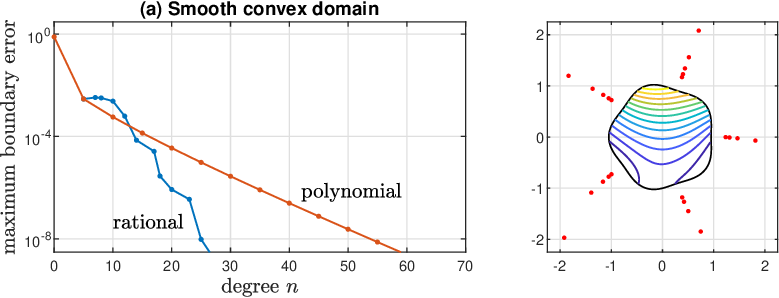}
\vskip 13pt
\includegraphics[scale=.85]{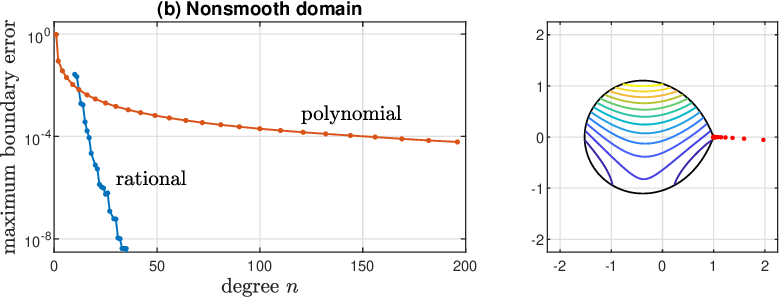}
\vskip 13pt
\includegraphics[scale=.85]{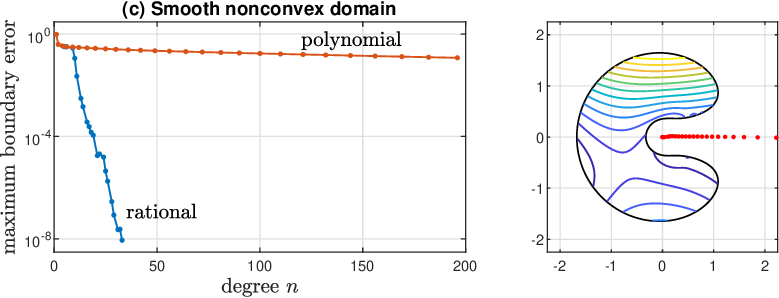}
\vskip 3pt
\end{center}
\caption{\label{figthree} Solutions of Laplace problems by
polynomials and rational functions of various degrees on three
domains.  On the smooth convex domain {\rm (a)}, both methods
converge rapidly.  On the domain {\rm (b)} with a corner singularity,
rational approximations are much more efficient than polynomials,
an effect going back to Newman in $1964$~{\rm \cite{newman}}.
On the smooth nonconvex domain {\em (c)} with an inlet, rational
approximations are again much more efficient, the new observation of
this paper.  Red dots mark the poles of rational approximations that
achieve error about $10^{-8}$.  These computations and comparisons
are made possible by approximation algorithms developed in the last
five years~{\rm \cite{VA,costa,aaa}}.} \end{figure}

Suppose we wish to solve numerically the Laplace problem
\begin{equation}
\Delta u = 0, ~z\in \Omega, \qquad
u = h, ~z\in \Gamma
\label{problem}
\end{equation}
in a simply-connected planar domain $\Omega$ bounded by an analytic Jordan curve $\Gamma$,
as suggested in Figure~\ref{figthree},
where $h$ is a real analytic function.\footnote{By
an analytic Jordan curve
we mean the one-to-one image of the unit circle under
an analytic function with nonvanishing derivative.  See~\cite[p.~2]{shapiro}
or~\cite[p.~2]{walsh}.}
(Everything can be generalized to less smooth geometries or data, to other types of
boundary conditions, and to domains with holes.  In Figure~\ref{figthree}(b), $\Gamma$
is not analytic.)
For convenience we think of $\Omega$ as complex, identifying $z = x + iy$.
An old idea, going
back to Walsh and Curtiss nearly a century ago~\cite{curtiss,walsh29},
is to approximate $u$ as the real part of a polynomial,
\begin{equation}
u(z) \approx \Re p(z),
\label{Rep}
\end{equation}
so that the problem reduces to the approximation of $h$ by $\Re p$ on $\Gamma$.
This idea builds on the facts that $u$ must be the real part of
a function $f$ that is analytic in $\Oc$~\cite{axler},
\begin{equation}
u(z) = \Re \kern -.5pt f(z),
\label{Ref}
\end{equation}
and that $f$ can be approximated on $\Oc$ by a polynomial,
\begin{equation}
f(z) \approx p(z).
\label{fp}
\end{equation}
In particular, Runge showed in 1885 that $f$ can be
approximated arbitrarily closely by polynomials in the supremum
norm~\cite{runge85}.

In the old days, computing good approximations $p$ would have been
problematic even if today's computers had been available, because
of the difficulty of finding well-distributed boundary points for
interpolation and also the lack of a well-conditioned basis in
which to represent the polynomial.  Today, the first difficulty
is bypassed by the use of least-squares fitting in a large number
of sample points on $\Gamma$, a proposal originating perhaps with
Moler in 1969~\cite{moler}, and the second is taken care of by the
use of the Vandermonde-with-Arnoldi method of on-the-fly Stieltjes
orthogonalization~\cite{VA}.  As a result, numerical solution of
planar Laplace problems by polynomial approximation is entirely
practical nowadays, so long as polynomials exist that approximate
$f$ efficiently.  For many problems this is an excellent numerical
method, quick and accurate and delivering a result that is exactly
analytic and trivial to differentiate.  Figure~\ref{figthree}(a)
is of this kind, and for further illustrations, see~\cite{series}
and~\cite[Figure~6.1]{VA}.

However, sometimes good approximating polynomials do not exist.
Figure~\ref{figthree}(b) shows a context in which this has been known
since the work of Bernstein, Jackson, and de la Vall\'ee Poussin
in the 1910\kern .3pt s, where $u$ has a boundary singularity.
Here the boundary contains a corner, and this prevents rapid
convergence by any polynomial approximations.

The first purpose of this paper is to show that the same effect
may arise even when $\Gamma$ and $h$ are analytic, as illustrated
in Figure~\ref{figthree}(c).  This is the situation when $\Omega$
is nonconvex, containing an inlet.  The figure shows that polynomial
approximations may converge at a negligible rate in such cases,
and we shall prove this mathematically.  Although the convergence is
exponential, the convergence constant is exponentially close to $1$,
rendering polynomial approximations useless in practice (Theorem~4).

On the other hand,
instead of approximating $u$ by a polynomial on $\Gamma$, one may approximate it by
a rational function with no poles in $\Omega$, 
\begin{equation}
u(z) \approx \Re r(z),
\label{Rer}
\end{equation}
which implicitly makes use of an approximation
\begin{equation}
f(z) \approx r(z).
\label{fr}
\end{equation}
Walsh considered this idea too~\cite{walsh29}, though it was even
further out of practical range in those days, with no robust
algorithms available for computing rational approximations,
even if suitable computers had been at hand.  Much more recently
Hochman, Leviatan, and White developed a method of this kind
in 2013~\cite{hlw}, and the AAA-least squares method, which
appeared a few years later, has made these computations quick and
easy~\cite{costa}, so that rational approximations are now a very
practical method for solving Laplace problems.

The second purpose of this paper is to prove that the speedups
possible with rational functions are transformative, not just in
cases like Figure~\ref{figthree}(b), where the power of rational
approximations has been known for a long time,
but also in cases like Figure~1(c), whose analysis is new.

Although this paper includes numerical illustrations to make
the points clear, its main purpose is theoretical.  The power of
rational approximations for solving Laplace and related problems
has been illustrated abundantly by numerical experiments in other
works~\cite{baddoo,costa2,costa,lightning,pnas,hlw,xue}.  The same
power applies to Helmholtz problems too~\cite{pnas}, although here
there is less literature and no theory.

One instance of a nonconvex domain is particularly tidy, the case
of an {\em inverted ellipse}, where everything can be worked out
explicitly.  For the inverted ellipse of parameter $\rho>1$, we
show that the degree of a polynomial approximation must be increased
asymptotically by $(\log(10)/1.16) \exp((\pi^2/4)/(\kern .7pt\rho-1))
\approx 2.0\cdot (11.8)^{1/(\kern .7pt\rho-1)}$ for each additional
digit of approximation accuracy (Theorem~5).  With $\rho=1.3$,
for example, hardly an extreme case, each digit of accuracy
requires an increase of the polynomial degree by about $7000$ (Table 1).
With rational approximations, on the other hand, each new digit
of accuracy requires an increase in the rational degree by just
$\sim {1\over 2}\log(10) /(\kern .7pt\rho-1) \approx 1.2/(\kern
.7pt\rho-1)$, so for $\rho = 1.3$, about $4$ rather than $7000$
(Theorem~9).  See Figure~\ref{fig3}.

There are a number of theoretical details in the upcoming pages,
but the essential argument can be summarized compactly, as follows.
Convergence rates of polynomial and rational approximations to
Laplace solutions depend on analytic continuation outside $\Gamma$.
This analytic continuation is described by the theory of the
Schwarz function, which asserts that there will usually be branch
point singularities not far from~$\Gamma$ (Theorems~1 and~2).
This is true for both convex and nonconvex domains, as one can
see by comparing the pole locations in Figures 1(a) and 1(c),
which in both cases are near the boundary but not extremely near.
The significance of nonconvexity arises at the next step of the
argument, where we track the consequences of these singularities.
For polynomial approximation, convergence is determined by a
conformal map of the exterior of\/ $\Gamma$ to the exterior of
the unit disk (Theorems 3--4 and Figures 3--4), and in the case of an inlet, the
singularity will map to a point exponentially close to the disk,
resulting in a convergence factor exponentially close to $1$
(section~\ref{inlets}).  For rational approximation, on the other hand,
poles can line up near branch cuts of the analytic continuation,
leading to a much faster exponential convergence rate associated
with a more forgiving doubly connected conformal mapping problem
onto an annulus (Theorems 6--8 and Figures 5--6).  In a word, the freedom of a
rational function to place poles near $\Gamma$ rather than just
at $\infty$ eliminates the ``crowding'' phenomenon that causes
exponential slow-down for polynomial approximations around inlets.

It is hard to find works in the literature that are close to the
present paper, particularly in combining ideas of analytic
continuation of Laplace solutions with polynomial and rational approximation
theory.  Two important related papers are by Millar~\cite{millar},
who applies the Schwarz function to track singularities of analytic
continuations of Helmholtz solutions
across boundaries, and Barnett and Betcke~\cite{bb}, who
apply the method of fundamental solutions to Helmholtz problems, also
with the aid of the Schwarz function, and show that good positioning
of the singularities outside the domain is crucial to obtaining
well-conditioned bases.  In these and other works related to analytic
continuation of solutions of elliptic PDE, it is common to formulate
the problem in terms of two independent complex variables, following
Vekua in the 1950\kern .3pt s~\cite{henrici,vekua}.  However, so far as I know,
this is not needed in the special case of the Laplace equation.

I must end this introduction with a historical and personal remark.
Most of the theory of this paper has connections with Joseph Walsh,
the great mid-20th century expert on polynomial and rational
approximation, who was a mathematics faculty member at Harvard during
1921--1966 until he retired and took a position at the University
of Maryland.  None of the algorithms that make these computations
practicable were available in his day, however, not to mention
computers of the necessary capabilities.  Today, 
we can apply this Walsh-style theory in
an entirely new environment.  Meanwhile my own career has brought
me to retirement after 26 years at Oxford, whereupon, in a strange
symmetry, I have taken a position at Harvard.  This is the first
paper I have written at my new institution.  I never met Walsh,
but I had indirect connections to him through his students Ted
Rivlin, Dick Varga, and Ed Saff, and with affection I think of the
present paper as a continuation of Walsh's research interests into
a computational era he could not have imagined.

\section{\label{schwarz}Analytic continuation and the Schwarz function}

The starting point of our analysis is analytic continuation across the boundary
curve $\Gamma$.  
Let $h$ be the boundary data function of $(\ref{problem})$, which we have assumed
is real analytic.  This implies that $h$ can be analytically continued to a complex
analytic function in a neighborhood of\/ $\Gamma$.
Similarly, let $f$ be the complex analytic function of (\ref{Ref}) whose real part is
the solution $u$ of the Laplace problem (\ref{problem}).  (The imaginary
part of $f$ is determined up to a constant, which plays no role in the
discussion.)  Like $h$, $f$ must extend to an analytic function
in a neighborhood of\/ $\Gamma$.
Along with $f$ and $h$, a third analytic function whose analytic
continuation across $\Gamma$ we shall work with is the difference,
\begin{equation}
d(z) = f(z) - h(z).
\label{diff}
\end{equation}
The important property of $d$ is that it is pure imaginary on
$\Gamma$, because $f$ and $h$ have equal real parts there.  Thus $h$ and 
$d$ are analytic functions in a neighborhood of $\Gamma$ that map
$\Gamma$ to subsets of the real and imaginary axes, respectively.  If\/ $\Gamma$ were a
straight line segment or an arc of a circle, we could now describe
the relationship of $h$ and $d$ inside and outside $\Gamma$ by the Schwarz
reflection principle.  In our more general case where $\Gamma$ is
an analytic arc, the generalization that comes into play is based
on what is known as the Schwarz function~\cite{davis,shapiro}.

The {\em Schwarz function} of\/ $\Gamma$ is defined as
the unique function $S(z)$ that is analytic in a neighborhood 
of $\kern .5pt\Gamma$ and takes the values $S(z) = \overline{z}$ 
for $z\in \Gamma$.
Note that $S(z)$ is determined by $\Gamma$, not $f$ or $h$ or $d$.
The significance of $S$ is that its complex conjugate, $\Sc(z)$, 
maps points close to $\Gamma$ on one side to their analytic reflections
on the other side.  Specifically, it is known that in a sufficiently
small neighborhood\/ $U$ of\/ $\Gamma$, 
$S$ is analytic and satisfies the reflection
property $\Sc(\Sc(z)) = z$.  Here is how $S$ enables analytic continuation:
if $h$ or $d$ is analytic in the part of\/ $U$ interior or exterior
to $\Gamma$, then it extends
analytically to the other side by the reflection formula
\begin{equation}
h(\kern .7pt\Sc(z)) =  \overline{h(z)}, \quad z\in U
\label{hcontin}
\end{equation}
or
\begin{equation}
d(\kern .7pt\Sc(z)) = - \overline{d(z)}, \quad z\in U,
\label{contin}
\end{equation}
respectively.
In words, the values of\/ $h$ and $d$ at $\Sc(z)$ are the reflections in the
real and imaginary axes, respectively, of their values at~$z$.  For details of these
developments see~\cite[chap.~6]{davis} or~\cite[Prop.~1.2]{shapiro}.

\begin{figure}
\vskip 8pt
\begin{center}
\includegraphics[scale=.80, trim=0 0 0 20, clip]{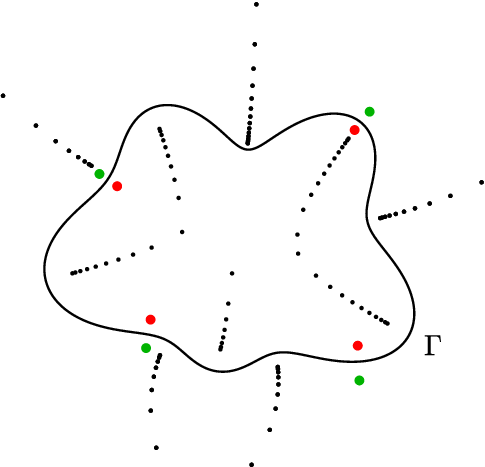}
\end{center}
\caption{\label{schwarzfig}An analytic Jordan curve $\Gamma$ and
the poles of a rational approximation on\/ $\Gamma$ to its Schwarz
function $S$.  The clustered strings of poles indicate that $S$ has
five branch points of $S$ near $\Gamma$ in the interior and another
five in the exterior.  The red and green dots mark four arbitrarily
chosen points on one side of\/ $\Gamma$ and their reflections on the
other side.  Analytic continuation of a function $d$ analytic and
imaginary on\/ $\Gamma$ is carried out by applying the reflection
condition $(\ref{contin})$ at such points.} \end{figure}

For an idea of the shape of the Schwarz function, consider
Figure~\ref{schwarzfig}.  This shows an analytic Jordan curve
$\Gamma$ for which a degree 107 AAA rational approximation $r$ to
$S$ has been computed.  (This took 2/3 s on our laptop based on the
Chebfun command \verb|[r,poles] = aaa(conj(Z),Z)|, where {\tt Z} is a vector
of 2000 sample points along the boundary.)  For $z\in\Gamma$, $r$
matches $S(z) = \overline{z}$ to accuracy $10^{-13}$.  The black dots
are the poles of $r$, and they give an indication of the behavior
of $S$ for $z\not\in \Gamma$.  Roughly speaking, with rational
approximations we expect poles to line up along the branch cuts of
the function being approximated.  More precisely, like all analytic
functions, $S$ does not intrinsically have branch cuts; these
are only introduced by humans when we wish to make the function
single-valued, or by rational approximations (in the approximate
fashion illustrated in the figure) as a by-product of their
near-optimality.  On the other hand $S$ has branch {\em points}, and
these are approximated by the cluster points of the poles of $r$.
In this case it appears that there are five branch points of $S$
near $\Gamma$ on the inside and five more on the outside.

Going out from $\Gamma$ beyond and around the branch points,
it would make sense to speak of a multi-valued analytic Schwarz
function, and this approach is taken in some of the theoretical
literature~\cite{shapiro} and will also turn up in our Theorem~8,
below.  However, multi-valued functions are mostly not relevant to
the present work, because polynomial and rational approximations are
single-valued.  Therefore, throughout this paper, our attention is
restricted to neighborhoods $U$ that are narrow enough or otherwise
confined in such a way that $S$ is analytic and single-valued.

Let $U$ be a neighborhood of\/ $\Gamma$ of this kind, with $S$
defined in $U$ satisfying $\Sc(U) = U$ and $\Sc(\Sc(z)) = z$, so
that $U$ is reflected into itself by $\Sc$.
We can now establish a foundational theorem for this paper
(compare~\cite[sec.~4]{millar}).

\medskip

{\em {\bf Theorem 1.~~Analytic continuation of \boldmath $f$ across
$\Gamma$.} Let $\Omega$, $\Gamma$, $f$, $h$, and $d=f-h$ be defined
as discussed above, and let\/ $U$ be a neighborhood of\/ $\Gamma$ in
which the Schwarz function $S$ of\/ $\Gamma$ is analytic and satisfies
$\Sc(U)=U$ and $\Sc(\Sc(z)) = z$.  If\/ $h$ is analytic in the part of\/ $U$ outside
$\Gamma$, then the same is true of $f$.}

\medskip
{\em Proof.}
Let $\Uin$ and $\Uout$ be the portions of\/ $U$ interior and exterior
to $\Gamma$, respectively.
If $h$ is analytic in $\Uout$, then by (\ref{hcontin}) it is analytic in $\Uin$.
By (\ref{diff}) it follows that $d$ is analytic in $\Uin$, since $f$ is
analytic throughout $\Omega$.
By (\ref{contin}), it follows that $d$ is analytic in $\Uout$.
By (\ref{diff}) again, it follows that $f$ is analytic in $\Uout$.
\qed

\medskip
To prove that polynomial approximation stagnates as illustrated in
Figure~\ref{figthree}, we will need a converse of this theorem.
We would like to show that when the Schwarz function $S$ has
singularities outside $\Omega$, they will usually block the analytic
continuation of the function $f$ associated with the Laplace problem
(\ref{problem}).  More precisely, I suspect it is true that if
$S$ cannot be analytically continued to a point $z_c$ outside of
$\Omega$, then the same will apply to $f$ if $h$
is a nonconstant function analytic outside $\Omega$.
(The situation is different if $h$ is a constant, since then $f$ will be a constant
too and thus analytically continuable to all of $\complex$, regardless
of the shape of $\Gamma$.)

It is worth mentioning that there are plenty of functions $h$ that are
real on~$\Gamma$ and analytic in $\complex\backslash\Omega$.
If\/ $\Gamma$ is the unit circle, then $x$
and $y$ are both functions in this class, since they can be written
$x = (z + z^{-1})/2$ and $y = (z - z^{-1})/2\kern .4pt i$.  The same therefore
applies on the unit circle to real polynomials in $x$ and $y$, and more generally, to
functions defined by Laurent series convergent for $0<|z| < \infty$
whose coefficients have the symmetry $a_{-k} = \overline{a_k}$.
For a general Jordan curve $\Gamma$, we may obtain an equivalent
class of functions $h(z)$ as transplants of these Laurent series
under a conformal map $\Phi$ of the exterior of $\Gamma$ to the
exterior of the unit disk.   (We will make use of $\Phi$ in the next
section in the context of polynomial approximation.)  Note that it
follows from here that, given any continuous real function $h_0$ on
$\Gamma$ and any $\varepsilon>0$, we can find a function
$h$ that is real on $\Gamma$ and analytic in $\complex\backslash\Omega$ with
$|h-h_0|<\varepsilon$ on $\Gamma$.

I don't know how to establish the singularity of $f$ in the
generality conjectured above.
Instead, following~\cite[sec.~4]{millar},
here is a statement restricted to the familiar situation in
which $S$ has a branch point as in Figure~\ref{schwarzfig} and in our later
example of the inverted ellipse.

\medskip

{\em {\bf Theorem 2.~~Limit to analytic continuation \boldmath
of\/ $f$ across \boldmath $\Gamma$.} Let\/ $\Omega$, $\Gamma$, $f$,
$h$, and $d$ be defined as discussed above, and let $z_c\in
\complex\backslash \Oc$ be such that $S$ can be analytically
continued up to $z_c$ but has a branch point there.
Then for some choices of\/ $h$
analytic in $\complex\backslash\Omega$, $f$ cannot be
analytically continued to a neighborhood of $z_c$.}

\medskip

{\em Proof.}
If $f$ can be analytically continued to a neighborhood
of $z_c$, then by (\ref{diff}) and the
analyticity of $h$, the same is true of $d$.
Let $U$ be a neighborhood of $\Omega$ as in
Theorem~1 that takes the form of an open set containing $z_c$ from which
has been removed
a branch cut $C$ extending from $z_c$ to the boundary.
Then $\Sc$ reflects $h$ and $d$ by (\ref{hcontin}) and
(\ref{contin}) from the portion of $U$ outside
$\Gamma$ to the portion inside.  
Since $z_c$ is a branch point,
$\Sc$ must take distinct values on the two sides of $C$.  Let $b\in C$
be such a point, mapping under $\Sc$ to two distinct points
$a_1,a_2\in \Omega$.
Since $d$ and $h$ are analytic at $b$, and both can be reflected
by $\Sc$ in $U$, we have $d(a_2) = d(a_1)$ and $h(a_2)=h(a_1)$, which implies
$f(a_2) = f(a_1)$ and $u(a_2) = u(a_1)$.  It follows from
the $|h-h_0|$ observation above, however, that this is not possible
for all choices of $h$.  We need only pick a boundary function $h_0$ for
which $u(a_1)\ne u(a_2)$, and then find an $h$ that is sufficiently
close to $h_0$.
\qed

\medskip

\section{\label{polytheory}Polynomial approximation}
Let $\En$ be the minimax error of degree $n$ polynomial approximation to
$f$ on $\Omega$,
\begin{equation}
\En = \inf_{p\in P_n} \|f-p\kern .7pt\|,
\label{Endef}
\end{equation}
where $\|\cdot\|$ is the supremum norm on $\Oc$ and $P_n$ denotes
the space of polynomials of degree $n$.
We now consider what Theorem 2 implies about the convergence rate of
$E_n$ to $0$ as $n\to\infty$.

\begin{figure}
\vskip 8pt
\begin{center}
\includegraphics[scale=.80, trim=0 0 0 0, clip]{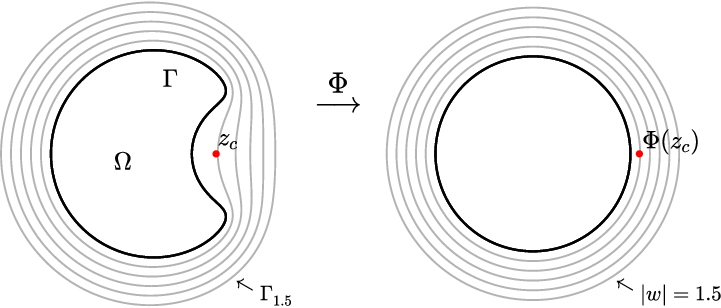}
\end{center}
\caption{\label{fig1}Illustration of the theory of polynomial approximation
of an analytic function $f$ on a
domain $\Omega$ bounded by an analytic Jordan curve\/ $\Gamma$.
The grey curves are level lines $\{\Gamma_r\}$ of the Green's function of\/ $\Gamma$,
$r = 1.1,1.2,\dots,1.5$,
which can be interpreted as preimages of circles outside the unit disk in a conformal
map $\Phi$ of the exterior of\/ $\Gamma$ to the exterior of the unit disk.  The analyticity
radius $R$ of\/ $f$ is the largest $r>1$ for which $f$ is analytic in the region
interior to $\Gamma_r$.  The red dot on the left marks a point on $\Gamma_R$ to which $f$ cannot
be analytically continued, and the absolute value of its image on the right is
the number $R$.   The convergence rate of best polynomial approximations is
given by $(\ref{rate})$.}
\end{figure}

Assume for simplicity that $f$ cannot be analytically continued
to all of $\complex$.  According to standard theory going back
to Walsh~\cite{walsh26,walsh} and presented beautifully by Levin
and Saff~\cite{ls}, the convergence rate is then exponential,
at a rate essentially $R^{-n}$ for some $R>1$.  The constant
$R$, which we shall call the {\em analyticity radius} of $f$ on
$\Omega$, has an interpretation in terms of a conformal map $\Phi$
of the exterior of\/ $\Omega$ to the exterior of the unit disk in
the $w$-plane with $\Phi(\infty) = \infty$.  For any $r\ge 1$, as
illustrated in Figure~\ref{fig1}, let $\Gamma_r$ be the preimage of
$|w|=r$ under $\Phi$.  Then $R$ is the largest $r$ for which $f$
is analytic in the region enclosed by $\Gamma_r$.\footnote{If
$z_c$ is a point on $\Gamma_R$ at which $f$ is not analytic,
then we can also write $R = \exp(\kern .7pt g(z_c))$, where $g$
is the Green's function of $\Oc$, that is, the harmonic function
defined in $\complex \backslash \Oc$ with $g(z) = 0$ on $\Gamma$
and with $g(z) \sim \log|z|$ as $z\to\infty$~\cite[Thm~3]{ls}.
The relationship between $g$ and $\Phi$ is $g(z) = \log |\Phi(z)|$.}

Here is Walsh's result.

\medskip

{\em {\bf Theorem 3.~~Polynomial approximation of an analytic function.}  Let the domain $\Omega$,
the function $f$ analytic in $\Omega$, the polynomial minimax error $E_n$, and
the analyticity radius $R$ of $f$ be defined as above.
The minimax approximation errors satisfy
\begin{equation}
\limsup_{n\to\infty} E_n^{1/n} = {1\over R}.
\label{rate}
\end{equation}
}

\medskip

We can paraphrase Theorem~3 like this:

\medskip

\begin{center}
\em\obeylines
The rate of polynomial approximation of\/ $f$ on\/ $\Omega$
is determined by its closest singularity outside\/ $\Omega$.
\end{center}

\medskip
\noindent ``Closest'' means in the sense of the level curves $\Gamma_r$, and the rate
in question is that of (\ref{rate}).  This is the ``Walsh half'' of the analysis of
our Laplace problem.
The ``Schwarz half'' was presented in Theorem~2 of the last section:

\medskip

\begin{center}
\em\obeylines
This closest singularity is normally the closest point\/ $z_c$ of nonanalyticity
of the Schwarz function $S$ (unless\/ $h$ has a singularity even closer).
\end{center}

\medskip
\noindent
Combining Theorems 2 and 3 gives us our basic result on slow
convergence of polynomial approximations of solutions of Laplace
problems.  As with Theorem~2, though this statement only mentions
``some'' choices of $h$, in practice it will be almost all of them.
The theorem references the analyticity radius $R$ of $S$, a notion
defined above for a function analytic throughout $\Omega$.  For a
function like $S$ analytic just in a neighborhood of $\Gamma$,
$R$ is the largest $r>1$ for which $S$ is analytic in the region
bounded between $\Gamma$ and $\Gamma_r$.

\medskip

{\em {\bf Theorem 4.~~Polynomial approximation of a Laplace solution.}  Let the Laplace problem
$(\ref{problem})$ with boundary data\/ $h$ have solution 
$u(z) = \Re \kern -.5pt f(z)$ as in $(\ref{Ref})$,
let\/ $R>1$ be the analyticity radius of the Schwarz function\/ $S$
on $\Omega$, and let\/ $E_n$ be the minimax error $(\ref{Endef})$
in degree $n$ polynomial approximation.
Then for some choices of\/ $h$ analytic throughout $\complex\backslash\Oc$,
and assuming $S$ has a branch point on the curve $\Gamma_R$,
these errors satisfy
\begin{equation}
\limsup_{n\to\infty} E_n^{1/n} = {1\over R}.
\label{rate2}
\end{equation}}

\medskip

\section{\label{invellpoly}Polynomials on the inverted ellipse}

\begin{figure}
\begin{center}
\vskip 14pt
\includegraphics[scale=.68, trim=0 0 0 0, clip]{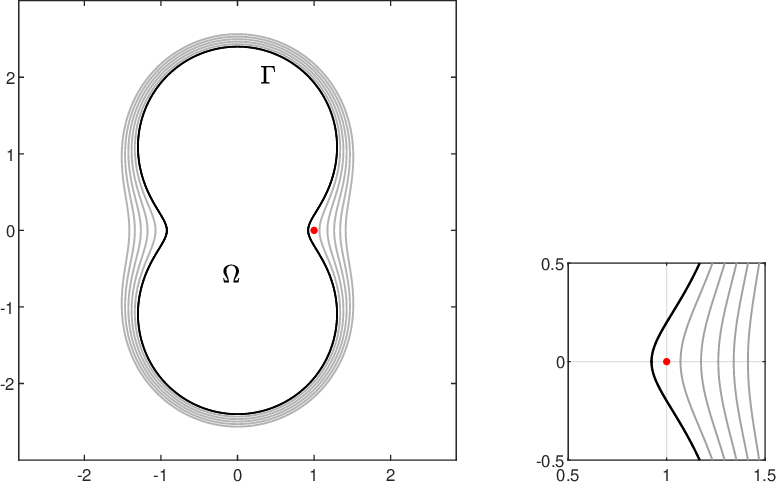}
\end{center}
\caption{\label{fig2}The inverted ellipse $I_\rho$ with parameter
$\rho = 1.5$, with level curves $\Gamma_r$ plotted
for $r = 1.02,1.04,\dots,1.10.$  The dot at $z=1$ corresponds
to the analyticity radius $r= R \approx 1.009$
at which the Schwarz function for\/ $\Omega$
has a branch point.  This value implies that although the boundary
of the region is analytic,
it will still take an increase
of the polynomial degree by about\/ $250$ asymptotically for
each additional digit of accuracy.  On the right, a close-up.}
\end{figure}

Theorem~4 has an elegant application to the special case in which
$\Omega$ is the region bounded by an inverted $\rho$-ellipse, as
sketched in Figure~\ref{fig2}.  For any $\rho>1$, we define the {\em
$\rho$-ellipse} $E_\rho$ to be the image in the $\zeta$-plane of the
circle $|w|=\rho$ in the $w$-plane under the Joukowsky map $J(w) =
(w+w^{-1})/2$.  Geometrically, $E_\rho$ is the ellipse with foci
$\pm 1$ whose semiminor and semimajor axis lengths sum to $\rho$.
The {\em inverted $\rho$-ellipse} is the reciprocal
$I_\rho = 1/E_\rho$, and our variables are related by $\zeta = z^{-1}$.

The Schwarz function for the inverted ellipse is known analytically.
For example on p.~25 of~\cite{davis} we find for the ellipse $E_\rho$
\begin{equation}
S(z) = \half(\kern 1pt\rho^2+\rho^{-2}) \kern .5pt z - \half (\kern 1pt\rho^2-\rho^{-2})\sqrt{z^2-1},
\label{SE}
\end{equation}
which implies that for $I_\rho$ we have
\begin{equation}
S(z) = \Bigl(
\half(\kern 1pt\rho^2+\rho^{-2}) \kern .5pt z^{-1} -
\half (\kern 1pt\rho^2-\rho^{-2})\sqrt{z^{-2}-1}\kern 2pt\Bigr)^{-1}.
\label{SI}
\end{equation}
For our purposes what
matters is that the only singularities of $S$ are a pair of branch points outside
$\Gamma$ at $z=\pm 1$ and a pair of simple poles inside $\Gamma$
at $\pm \kern .4pt i\kern .4pt (\kern 1pt \rho^2 - \rho^{-2}\kern .7pt )/2$.
If branch cuts are drawn along $[1,\infty)$ and $(-\infty,-1]$, then $S$ becomes 
meromorphic in the remaining slit domain
$\complex \backslash \{ [1,\infty)\cup (-\infty,-1]\}$, analytic in the portion
outside $\Omega$.

To determine the analyticity radius for $I_\rho$, we accordingly need to know
the image of $z=1$ under a conformal map of the exterior of $I_\rho$
to the exterior of the unit disk.
The conformal map of the interior of $I_\rho$ to the unit disk is 
elementary: it reduces to the map of the exterior of $E_\rho$, which is
essentially the Joukowsky map.  The conformal map
of the exterior of $I_\rho$, on the other hand, reduces to the map of the
interior of $E_\rho$, which is not elementary.
The required formula involving a Jacobian
elliptic function was derived
by Schwarz in 1869 and is presented in a number of sources
including~\cite{nehari}, \cite{szeg},
and Figure~3.2 of~\cite{ht}.
According to equation (21) of~\cite{szeg}, the focus $\zeta=1$ of the ellipse 
maps to the point
\begin{equation}
\varphi^{-1}(1) = { \rho^{-1^2} + \rho^{-3^2} + \rho^{-5^2} + \cdots \over
 {1\over 2} + \rho^{-2^2} + \rho^{-4^2} + \rho^{-6^2} + \cdots}.
\end{equation}
The analyticity radius we need follows by taking the reciprocal.

\medskip

{\em {\bf Theorem 5.~~Polynomial approximation on the
inverted ellipse.}  If\/ $\Omega$ is the domain bounded by the inverted
$\rho$-ellipse $I_\rho$ for some $\rho>1$, the associated analyticity
radius $R$ of Theorem~$4$ is
\begin{equation}
R = {{1\over 2} + \rho^{-2^2} + \rho^{-4^2} + \rho^{-6^2} + \cdots \over
\rho^{-1^2} + \rho^{-3^2} + \rho^{-5^2} + \cdots }.
\label{Rparam}
\end{equation}
Asymptotically as $\rho\to 1$, $R$ satisfies
\begin{equation}
R-1 \sim 4 \kern .4pt e^{-\pi^2/(4\log \rho)} 
\sim A \kern .4pt e^{-\pi^2/(4(\kern .7pt\rho-1))} 
\label{ierate}
\end{equation}
with $A= 4\kern .7pt \exp(-\pi^2/8)\approx 1.16485$.
}

\medskip

{\em Proof.} The derivation of (\ref{Rparam}) was given in the discussion
above, and I am grateful to Jon Chapman of Oxford and Alex Barnett of the Flatiron
Institute for proofs of (\ref{ierate}).  The following particularly elegant
argument comes from Barnett.  If we double both the numerator and the
denominator of (\ref{Rparam}) by extending the sums to $-\infty$,
the quotient is unchanged in value and takes the simple form
\begin{displaymath}
R = \sum_{k=-\infty}^\infty \rho^{-4k^2} \kern -4pt \left / 
\kern -3pt \sum_{k=-\infty}^\infty \rho^{-4(k+{1\over 2})^2} \right.,
\end{displaymath}
or equivalently, on setting $a = 4\log \rho$,
\begin{equation}
R = {A\over B} = \sum_{k=-\infty}^\infty \exp(-ak^2)\kern -2pt \left / 
\kern -3pt \sum_{k=-\infty}^\infty \exp(-a(k+\textstyle{1\over 2})^2)\right..
\end{equation}
This is the ratio of two infinite trapezoidal quadrature approximations
to the Gaussian, which explains
why it converges exponentially to 1 as $a\to 0$.
By the Poisson summation formula applied to the
Fourier transform pair $f(x) =\exp(-ax^2)$ and
$\hat f(k) =c \exp(-\pi^2 k^2/a)$,
where $c$ is a constant whose value does not matter, we have
for a new constant $c'$
\begin{displaymath}
A =c'\left(1 + 2\kern .5pt e^{-\pi^2/a} + 2\kern .5pt e^{-(2\pi)^2/a} + \cdots\right),
\end{displaymath}
and, since the translation by $1/2$ leads to alternating signs in the Poisson
summation formula,
\begin{displaymath}
B =c'\left(1 - 2\kern .5pt e^{-\pi^2/a} + 2\kern .5pt e^{-(2\pi)^2/a} - \cdots\right).
\end{displaymath}
These formulas imply $A/B - 1 \sim 4e^{-\pi^2/a}$ as $a\to 0$.
\qed

\medskip

The exponential dependence of $R$ on $\log \rho$ or $\rho-1$ in Theorem~5
is striking.  Table~1 lists $R$ for various
values of $\rho$ decreasing toward $1$.   The final column shows
$\log(10)/\log(R)$, the increase of degree required asymptotically
for each improvement of accuracy by one digit.

\vskip 8pt 

\begin{table}[h]
\begin{center}
\caption{\label{table1}Analyticity radii and corresponding 
convergence rates for polynomial
approximations to solutions of Laplace problems on the inverted ellipse $I_\rho$ for
various $\rho$.  All numbers are rounded to $2$ significant figures.}
{\small
\begin{tabular}{c l r}
$\rho$ & $~~~R$ & degree increase per digit \\
\hline\\[-7pt]
$ 2 $ & $1.12            $ & $20$ \\
$1.9$ & $1.089           $ & $27$ \\
$1.8$ & $1.062           $ & $38$ \\
$1.7$ & $1.038           $ & $60$ \\
$1.6$ & $1.021           $ & $110$ \\
$1.5$ & $1.0091          $ & $250$ \\
$1.4$ & $1.0026          $ & $880$ \\
$1.3$ & $1.00033         $ & $7000$ \\
$1.2$ & $1.0000053       $ & $430{,}000$ \\
$1.1$ & $1.000000000023  $ & $100{,}000{,}000{,}000$
\end{tabular}
}
\end{center}
\end{table}

\vskip 8pt 

\section{\label{rattheory}Rational approximation}
Polynomials are rational functions all of whose poles are constrained
to lie at $z=\infty$.  When rational functions without this
constraint are allowed, the approximation power increases enormously
for the kinds of domains we are interested in.  This happens because
poles can now line up along branch cuts of $f$, or, equivalently in
most cases, along branch cuts of the Schwarz function $S$.  We saw
this effect in Figure~\ref{figthree}.  The phenomenon of poles
approximating branch cuts is well known in rational approximation
theory, associated in particular with theorems of Stahl for Pad\'e
approximation~\cite{stahl}, and is illustrated, for example, in the
papers~\cite{baddoo,costa2,costa,lightning,pnas,hlw,aaa,analcont,xue}.

An explanation of this behavior and of the approximation power of
rational functions comes again from Walsh~\cite{walsh}.  In Walsh's
theory, key roles are played by the $n$ poles $\{\pi_k\}$ of a
rational function $r$ of degree $n$, which can be thought of as
point charges for a two-dimensional potential $\log|z-\pi_k|$, and
any $n$ interpolation points $\{z_k\}$ where $f(z_k)=r(z_k)$, which
can be thought of as point charges of opposite sign associated with
potentials $-\log |z-z_k|$.  By considering these potentials in the
context of the Hermite integral formula representing $f(z)-r(z)$
at other points $z\in\complex$, one can derive error bounds
associated with potential theory.  A few decades after Walsh,
Gonchar and others took the theory further, showing how in the
limit $n\to \infty$, minimax approximation rates are governed by
the equilibrium potentials for a continuum distribution.

We will not go into the details of this theory, which are summarized
in~\cite[pp.~91--93]{ls} and more fully in sections 6--8 of~\cite{analcont}.
We will, however, state its main conclusion, which is illustrated in
Figure~\ref{ratfig}.  Let $f$ be analytic on $\Omega$ and
extend analytically (i.e., as a single-valued analytic function) to some larger domain 
$\Omega'\subseteq \complex$ enclosing $\Oc$ in its interior.  For simplicity we suppose
that $\Omega'$ is also a Jordan domain, with boundary $\Gamma'$.  Then $\Gamma$ and
$\Gamma'$ are the inner and outer boundaries of an annular region $K$, called
a {\em condenser}.  Such a region can be conformally mapped onto a circular annular
region in $\complex$ whose inner boundary
is the unit circle and whose outer boundary is
the circle about the origin of some radius $R^*>1$,
the {\em modulus} of the condenser.  The number $R^*$ is uniquely determined, and the conformal
map itself is also uniquely determined up to a rotation.  As in 
(\ref{Endef}), but now for rational approximations, we define
\begin{equation}
\Enn = \inf_{r\in R_n} \|f-r_n\|,
\label{Endefr}
\end{equation}
where $\|\cdot\|$ is the supremum norm on $\Oc$ and $R_n$ denotes the
space of rational functions of degree $n$ (i.e., with at most $n$
poles counted with multiplicity, including any poles at $\infty$).
Here is Walsh's result:

\begin{figure}
\vskip 8pt
\begin{center}
\includegraphics[scale=.77, trim=0 0 0 0, clip]{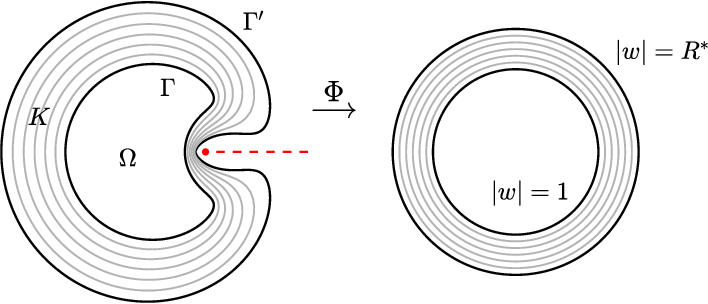}
\end{center}
\caption{\label{ratfig}A rational companion to Figure\/ $\ref{fig1}$,
illustrating the theory of rational approximation of an analytic function $f$ on a
domain $\Omega$ bounded by an analytic Jordan curve\/ $\Gamma$.  We suppose $f$ has
a singularity outside $\Gamma$, marked by the red dot, but that $f$ can be analytically
continued to a condenser region $K$ bounded between $\Gamma$ and a larger Jordan
curve $\Gamma'$ that avoids this singularity.  If $\Phi$ is a conformal map of the
condenser onto a circular annulus $1<|w|<R^*$, then rational approximations
to $f$ can converge at the exponential rates {\rm (\ref{rrate})} or
{\rm (\ref{rate5})} determined by $R^*$.  They
achieve this by placing poles approximately along branch cuts, as suggested
by the dashed line.  Compare Figure {\rm \ref{figthree}(c)}.  The conformal
map for this image was computed by the AAA-least squares method.}
\end{figure}

\medskip

{\em {\bf Theorem 6.~~Rational approximation of an analytic function.}  Let the domain
$\Omega$, the function $f$, the rational minimax error $\Enn$, the
condenser $K$, and the
condenser modulus $R^*$ be defined as above.
The minimax errors satisfy
\begin{equation}
\limsup_{n\to\infty} \Enn^{1/n} \le {1\over R^*}.
\label{rrate}
\end{equation}
}

\medskip

Recall that Theorem 3 for polynomials involved a number $R$ that
was fully determined by $f$ and $\Omega$, defined as the largest
$r>1$ for which $f$ could be extended analytically to within the
Green function contour $\Gamma_r$.  Here with rational functions,
by contrast, the number $R^*$ as we have defined it depends on the
choice of the condenser $K$ to which $f$ is analytically
continued.  In cases where $f$ has branch points, a single-valued
continuation will need to be restricted to a domain that avoids
these, and some choices will lead to larger values of $R^*$ than
others.  An optimal choice can be made, leading to a maximal value
of $R^*$, but we will not go into that here, simply accepting that
Theorem~6 holds for any choice of $K$ and its associated $R^*$.
As a consequence the next two theorems involve inequalities, not equalities.

The crucial point is that $R^*$ can be much bigger than $R$, allowing
rational approximations to converge much faster than polynomials.
In particular we will see in sections~\ref{invellrat} and~\ref{inlets}
that this happens
with domains $\Omega$ involving inlets, where $R$ is exponentially
close to $1$ whereas $R^*$ may be only algebraically close.

\medskip

{\em {\bf Theorem 7.~~Rational approximation of a Laplace solution.}  Let the Laplace problem
$(\ref{problem})$ with boundary data\/ $h$ have solution 
$u(z) = \Re \kern -.5pt f(z)$ as in $(\ref{Ref})$.
Let the Schwarz function\/ $S$ of\/ $\Gamma$ be
analytically continuable to a condenser\/ $K$ about\/ $\Omega$, as
defined above, with modulus $R^*$.  Assume that $\Sc$ reflects $K$
into a subset of\/ $\Omega$ such that $U = K\cup \Gamma \cup \Sc(K)$,
as in Theorem\/~$1$, is a neighborhood of\/ $\Gamma$ in which $S$ is
analytic and satisfies $\Sc(U) = U$ and $\Sc(\Sc(z)) = z$.
If $h$ is analytic in $K$, then
the rational minimax errors $\Enn$ of $(\ref{Endefr})$ satisfy
\begin{equation}
\limsup_{n\to\infty} E_{nn}^{1/n} \le {1\over R^*}.
\label{rate3}
\end{equation}
}

\medskip

{\em Proof.}  By Theorem~1, $f$ is analytic in $K$, and (\ref{rate3}) follows
from Theorem~6. \qed

\medskip
The $1/R^*$ convergence factor of Theorems 6 and 7 is pessimistic
in many cases: the actual convergence factor is often $1/(R^*)^2$.
Equivalently, we may say that $1/R^*$ is often an upper bound
not just for $\limsup_{n\to\infty} \Enn^{1/n}$ but also for
$\limsup_{n\to\infty} \Enn^{1/2n}$.  Intuitively speaking, this
factor of 2 is associated with the fact that rational functions have
twice as many parameters as polynomials, a property not exploited in
the theory leading to Theorems 6 and 7.  These issues are summarized
in section 7.5 of~\cite{analcont}, whose substance is derived from
the important paper~\cite{rakh} by Rakhmanov.  According to what
Rakhmanov calls the {\em Gonchar--Stahl $\rho^2$-theorem,} the factor
of 2 speedup always applies in the case of a function that can be
analytically continued to a multivalued analytic function along
all curves in\/ $\complex$ that avoid a certain fixed set $\Sigma\subseteq \complex$
of capacity zero.  In particular, it applies when $f$ is analytic
apart from branch points.  We summarize the improvement to Theorems~6
and~7 as follows.  \medskip

{\em {\bf Theorem 8.~~Factor of 2 speedup for functions with algebraic
branch points.}  The convergence rates of
Theorems $6$ and $7$ can be doubled to 
\begin{equation}
\limsup_{n\to\infty} E_{nn}^{1/n} \le {1\over (R^*)^2}.
\label{rate5}
\end{equation}
if\/ $f$ can be analytically continued to a multivalued analytic
function along all curves in $\complex$ that avoid a fixed
set\/ $\Sigma\subseteq\complex$ of capacity zero.}

\medskip

{\em Proof.}  The estimate (\ref{rate5}), which is the upper-bound half
of the Gonchar--Stahl theorem, originates in~\cite{stahl86}.  \qed

\medskip

\section{\label{invellrat}Rational functions on the inverted ellipse}
\begin{figure}[p]
\begin{center}
\vskip 8pt
\includegraphics[scale=.88]{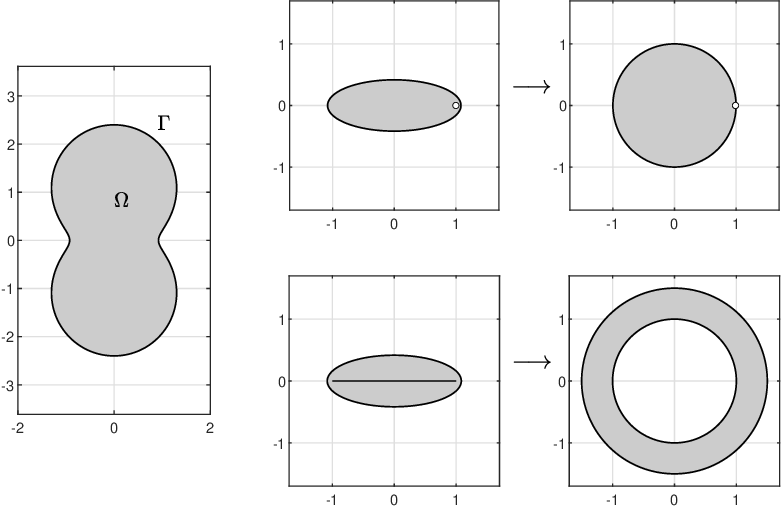}
\end{center}
\caption{\label{fig0} A summary in images of polynomial and rational
convergence rates for the inverted ellipse $I_\rho$ as discussed in
sections~$4$ and\/~$6$.
On the left, the domain\/ $\Omega$ bounded by $I_\rho$ with
$\rho = 1.5$.  On the right, the conformal maps
of two domains bounded by the corresponding ellipse $E_\rho$,
first without and then with 
a slit along the midline $[-1,1]$.  According to the arguments
leading to Theorems~$5$ and\/~$9$ and summarized in
Figures\/~$\ref{fig1}$ and\/~$\ref{ratfig}$, these maps determine
the exponential convergence rates for polynomial and rational
approximation methods, respectively,
for solving Laplace problems on $\Omega$.  The polynomial rate
is very slow because the image of $z=1$ in the upper map is exponentially close to $1$
(about $0.9909$ for this choice of $\rho$, marked by a white dot).
The rational rate is much faster because the outer boundary of the annulus in the
lower map is only algebraically close to the inner one (distance $0.5$).}
\end{figure}

We now apply Theorems 7--8 to the inverted ellipse $I_\rho$.  This proves
surprisingly easy.  We saw in section 4 that the conformal map
of the exterior of $I_\rho$ to the exterior of the unit disk is not
elementary.  For rational approximation, however, since the Schwarz
function is analytic in the slit domain
$\complex \backslash \{ \Oc \cup [1,\infty)\cup (-\infty,-1]\}$, the conformal map
we need is of the exterior of $\Oc$ in this slit domain onto a circular annulus.  
Taking reciprocals, what is at issue is the conformal map onto a circular annulus of
the doubly-connected domain bounded by the $\rho$-ellipse $E_\rho$ and slit along $[-1,1]$.
This is essentially just the Joukowsky map, as summarized in Figure~\ref{fig0}.

The gives us the following theorem.

\medskip

{\em {\bf Theorem 9.~~Rational approximation on the
inverted ellipse.}  If\/ $\Omega$ is the domain bounded by the inverted
$\rho$-ellipse $I_\rho$ for some $\rho>1$, then if\/ $h$ can be analytically continued
to all of\/ $\complex\backslash\Omega$, the associated modulus $R^*$ of
Theorem~$8$ can be taken to be any value $R^*<\rho$.
It follows that for such $h$, minimax rational approximation errors for
solutions of the Laplace problem satisfy
\begin{equation}
\limsup_{n\to\infty} E_{nn}^{1/n} \le {1\over \rho^2}.
\label{rate4}
\end{equation}
}

\medskip

\begin{figure}[p]
\vskip 10pt
\begin{center}
\includegraphics[scale=.45]{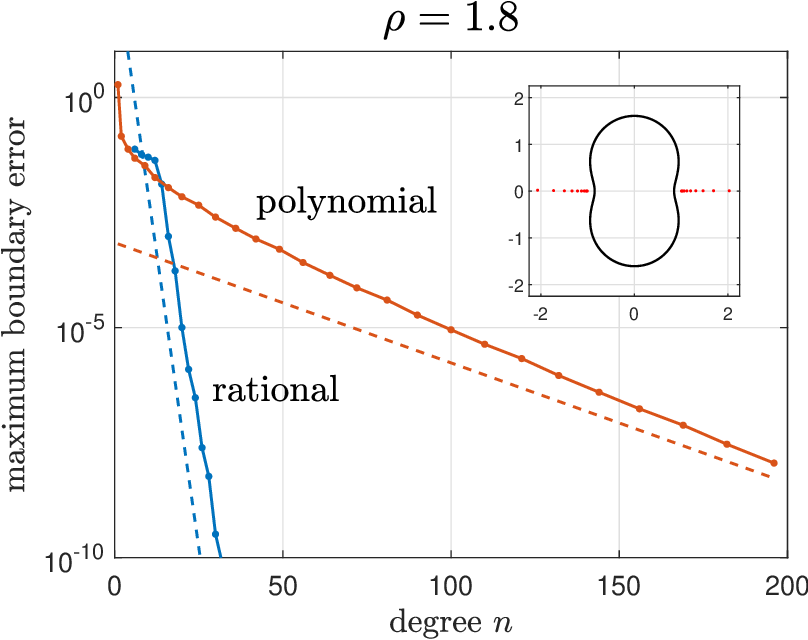}~~~~~~\includegraphics[scale=.45]{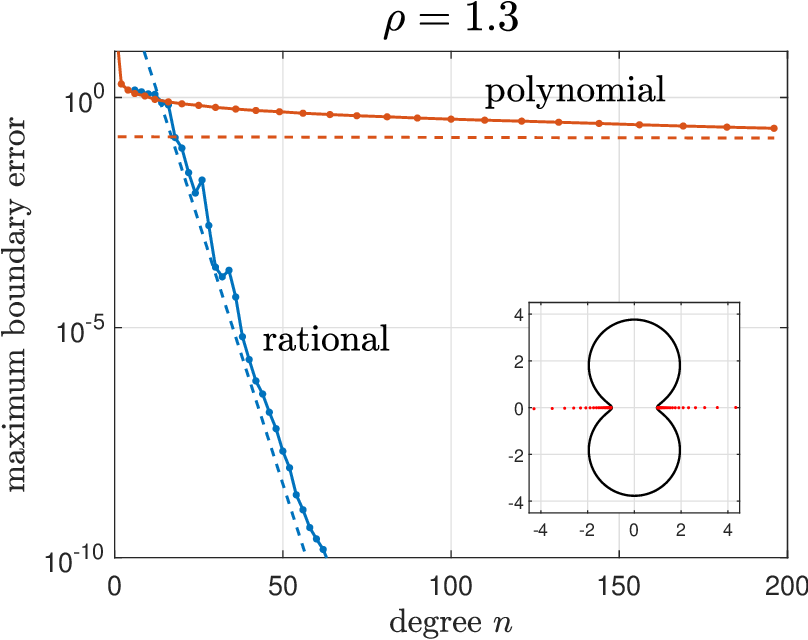}
\end{center}
\caption{\label{fig3} Numerical confirmation of Theorems $5$ and\/ $9$
for solving a Laplace problem by polynomial and rational
approximation on the domain bounded by the inverted $\rho$-ellipse $I_\rho$ (Vandermonde
with Arnoldi least squares~{\rm \cite{VA}} and AAA-least
squares\/~{\rm\cite{costa}}, respectively).
With $\rho = 1.8$, the indentation is so mild that exponential convergence
is observed for both methods, though at a lower rate for
polynomials than rational functions.
With $\rho = 1.3$ the exponential convergence rate for polynomial approximation is
essentially zero and we see subexponential convergence, as if\/ $\Omega$
had a corner.  The rational method still converges rapidly, and the poles marked
as red dots in the inset show that this is achieved by means of poles
delineating approximate branch cuts $(-\infty -1]$ and $[1,\infty)$.
Dashed lines show the theoretically predicted slopes
$(\ref{ierate})$ and\/ $(\ref{rate4})$; their heights are arbitrary.}
\end{figure}

Figure~\ref{fig3} shows a pair of numerical experiments confirming Theorems~5
and~9.  The figure has two panels,
corresponding to domains $\Omega$ bounded by inverted ellipses
$I_\rho$ with $\rho = 1.8$ and $\rho = 1.3$.  Each panel compares the
convergence of a polynomial approximation method (Vandermonde with
Arnoldi least squares~\cite{VA}) against a rational one
(AAA-least squares~\cite{costa}).
The boundary data function is $h(x+iy) = (y+1)^2$, with the shift by
$1$ introduced to break symmetry, and the boundary is discretized
by $1500$ points, so that these computations involve matrices of reasonably
modest size with $1500$ rows and a few hundred columns.

For $\rho=1.8$, where $\Omega$ is barely indented at all, Figure~\ref{fig3} shows that the
polynomial method converges about $5$ times more slowly than the
rational method, requiring an increase of $n$ by about 38 for each
digit of accuracy as listed in Table~1.

The case $\rho=1.3$ shows the extreme behavior predicted by Theorem~5.  
The rational method has slowed down by a factor of 2, whereas
the polynomial method has slowed down by a factor of 200, 
now needing an increase of degree by $7000$ for each additional digit of
accuracy.  Note that convergence of the polynomial method is 
still observed, but showing the upward-curving form associated with
subexponential convergence for a non-smooth boundary rather
than a straight line for exponential behavior.  Evidently it
would be impossible in practice to get even two digits of accuracy by this method,
even though the boundary is analytic and the Schwarz function extends 
analytically a nonnegligible distance $0.0346$ outside it.  This distance is just
algebraically small, but its consequences are exponential.

The computations for Figure~\ref{fig3}, about 120 numerical Laplace solutions all together,
required 9 seconds on our laptop.

\section{\label{inlets}Approximation on a general domain with an inlet}

The exponential effect that makes polynomial approximation on nonconvex
domains problematic goes by the name of ``crowding'' in the literature
of numerical conformal mapping.
(In elasticity theory it is called Saint-Venant's principle.)
The observation here is that a conformal
map involving a long and narrow peninsula or finger,
which hardly needs to be very
long and narrow, will involve exponential distortions.  Peninsulas become
inlets in our context because approximation on a domain $\Omega$ depends on
the conformal map of its complement $\complex\backslash\Oc$.

The most systematic analysis of crowding that I know of appears
in~\cite{cm1}.  The following definition is used there for
an analytic Jordan domain $\Omega$:
\medskip

{\narrower
We say that $\Omega$ {\em contains a finger of length
$L>0$} if there is a rectangular channel
of width $1$ defined by a pair of parallel line segments of
length $L$, disjoint from $\Omega$, such that $\Omega$ extends
all the way through the channel with parts of $\Omega$ lying outside both
ends.  
\par}

\medskip

\noindent On the basis of this definition it is shown that the associated
harmonic measure (Theorem 2 of~\cite{cm1}), conformal mapping derivative (Theorem 3), 
radius of univalence (Theorem 4), and approximating
polynomial degrees (Theorem 5) all scale in ways controlled by the
factor $\exp(\pi L)$.  Thus, for example, peninsulas of
length-to-width ratios $1$, $2$, and $3$ induce distortions of magnitudes
on the order of $23$, $540$, and $12{,}000$.

In the context of the present paper, it follows from these theorems
that the degrees of polynomial approximations to solutions of
Laplace problems on any domain $\Omega$ will have to grow by a factor at least
as large as order $\exp(\pi L)$ for each additional digit of accuracy, if $\Omega$
contains an inlet of length-to-width ratio $L$.  We do not
spell out precise theorems.

It is interesting to check how closely this general result matches
Theorem~5 for the case of the inverted ellipse $I_\rho$. 
Setting $\varepsilon = \rho -1$,
we find that for $I_\rho$ the inlet parameter $L$ scales as\footnote{Equation
(\ref{fingerlength}) is derived by considering the rectangle
$1\le \Re z \le 2$, $-4\kern .4pt \varepsilon \le \Im z \le 4\kern .4pt \varepsilon$, with
aspect ratio $L = 1/8\kern .4pt \varepsilon$.  For $\varepsilon < 0.225$,
the right inward-pointing finger of $I_\rho$ is
contained within the top and bottom sides of this rectangle while extending outside both the
left and right ends.}
\begin{equation}
L \sim {1\over 8\kern .4pt \varepsilon}, \quad \varepsilon \to 0.
\label{fingerlength}
\end{equation}
This implies by the $\exp(\pi L)$ results above that
the analyticity radius $R$ of sections~\ref{polytheory} and~\ref{invellpoly} satisfies
\begin{equation}
R - 1 \le O(\exp(-\pi/8\kern .4pt \varepsilon)),
\end{equation}
which is looser by the factor $2 \pi$ than the actual result of (\ref{ierate}),
\begin{equation}
R - 1 \sim O(\exp(-\pi^2/4\kern .4pt \varepsilon)).
\end{equation}

We have seen that polynomial approximations of Laplace solutions can be expected
almost always to be ineffective on domains with inlets.  Conversely,
will rational approximations almost always be effective?
Certainly not always, for the boundary may be nonsmooth in the sense of
having singularities of the Schwarz function very close by.  But as a smooth example
it is interesting to consider the inverted ellipse once more.
From (\ref{fingerlength}), the finger length is about
$1/8\kern .4pt \varepsilon$, whereas from
Theorem~9, the convergence constant is $\rho^2 \sim 1+ 2\kern .4pt \varepsilon$.
Thus in this case at least, the convergence rate for rational
approximations slows down only linearly with the length of the finger.

\section{\label{dicuss}Discussion}
We have shown that polynomial approximations to solutions of Laplace problems
are essentially useless on non-convex domains, whereas rational approximations
are often very effective.  In principle these results might
have been obtained decades ago, but in practice, everything has changed since
the AAA and AAA-least squares algorithms have made rational functions an easy
tool for numerical computation.

The Laplace equation $\Delta u = 0$ is the limit
of the Helmholtz equation $\Delta u + k^2 u = 0$ as $k\to 0$, and poles
of rational functions can be regarded as the
limits of point singularities of certain Hankel functions~\cite{pnas}.   It is
known to experts
in the solution of Helmholtz equation by the method of fundamental solutions that when
a region has an inlet, it is crucial to include some Hankel
singular points therein~\cite{bb}.
Thus the present paper could be regarded as investigating a limiting case of
an effect that is known in the Helmholtz context, though we are not aware
that theorems are to be found in that literature analogous to what we have
developed here.

\section*{Acknowledgments} 
I am grateful for helpful comments (and in two cases for the derivation
of equation (\ref{ierate}) in 
Theorem~5) from Alex Barnett, Jon Chapman, Peter Ebenfelt, Astrid
Herremans, Daan Huybrechs, Stefan Llewellyn Smith, Ed Saff, and Yidan Xue.

\end{document}